\documentclass{notices}

\title{
Mathematics and Dance:
\\ Notes from an emerging interaction
}

\author{
  Reggie Wilson
  \affil{
    The first author is the Artistic Director and Choreographer at the Fist and Heel Performance Group. His email address is rwilson@fistandheelperformancegroup.org}
  \and
  Jesse Wolfson
  \affil{
    The second author is an Associate Professor of Mathematics at the University of California, Irvine.  His email address is wolfson@uci.edu.
   }
}

\begin{document}

\maketitle

What can math do for dance? What can dance do for math? For almost 10 years, we—Choreographer Reggie Wilson and Mathematician Jesse Wolfson—have been exploring and investigating these questions, both in the rehearsal studio, the theater and on zoom with my (Wilson’s) \href{http://www.fistandheelperformancegroup.org/}{Fist and Heel Performance Group}, and also in the classroom and university research community with my (Wolfson’s) students and colleagues. We aim to share with you here some of the answers that are beginning to emerge from our exploration.

It can help to start with definitions.  William Thurston \cite{Th} gives a rough definition of {\em math} as ``the theory of formal patterns.'' In a different voice, we might say that math is a sustained exploration of our intuitions of space and number and of how we make sense of them.  By {\em dance}, we mean the motions, actions, movements, (and kinesthetic experience) of our own and other bodies in 3-dimensional space, organized intentionally over time, improvisational or predetermined (set).   

Dance and mathematics represent two of the oldest practices through which humans access and engage with space and how bodies move through it.  In our culture and the lineages we often trace for it, these two practices have remained almost entirely distinct, but, a survey of contemporary activity reveals multiple collaborations and cross-currents linking dance and math.  We believe these interactions signify a potentially transformative shift.  As in the early phases of any new technology, we don’t yet understand the full scope of how the interaction between dance and mathematical sciences can be valuable. However, we can make initial assumptions based what has arisen, including applications to: K-12 education, concert dance, robotics, mathematical disciplines such as dynamics, topology and geometry, and to creating and strengthening linkages—across class, culture, and geography—building political and social agency.  We hope this article can help bring awareness to these interactions and inspire others.

\section*{Examples and Non-Examples}
Some examples and non-examples might help to illustrate what we mean. Let’s start with the latter.  Dance can be a medium to communicate mathematical ideas, as in  \href{https://nancyscherich.com/}{Nancy Sherich}'s (Toronto, Dept. of Math, Postdoctoral Scholar) \href{https://www.youtube.com/watch?v=MASNukczu5A&feature=youtu.be}{winning entry in the 2017 “Dance your PhD” contest}.  Mathematics can be a muse for generating new dance, as in the work of Merce Cunningham and John Cage.  Or, math and dance can be sources of mutual analogies and resonances that inform how one perceives and creates work in each separately, as described in Sophia Merow’s 2019 Notices article e \href{https://www.ams.org/journals/notices/201902/rnoti-p247.pdf?trk=1802&cat=collection}{``Making dance, making math: parallels.''}  These are interesting, useful and valid as answers to the questions we started with.  They are also not what we want to focus on today.  

What we have in mind are encounters between mathematics and dance as bodies of knowledge and practice, where the knowledge and practice of each can be used to solve (and pose) problems of intrinsic interest for the other. Here are some examples drawn from our experience.  This list is not exhaustive!

\paragraph*{I. Fractals.} In 2012, things got rolling as I (Wilson) and Fist and Heel were in the process of sourcing and developing movement for a new performance work.  Many of the movements and performance practices I was engaging with exhibited Africanist formal features—rhythms, brief sequences of movement, striking variations in quality and force or tempo—that did not map easily onto the schema and formalisms that characterize traditional Western concert dance.  Around this same time, I encountered Ron Eglash’s {\em African Fractals} \cite{Eg} which documents extensive use of fractal symmetry in African material cultures. I hypothesized that fractals were also prevalent in the performative cultures of Africa, and I invited Wolfson to help ``translate'' Eglash’s text and to give his perspective on a selection of Africanist forms of music and dance that had most piqued my curiosity to possible usage of things fractal. These consultations were foundational for the successful creation of my evening-length contemporary dance work, {\em Moses(es)}, and they laid the groundwork for our ongoing conversations and thinking about the relationships of dance and math.  These engagements were singularly focused on exploring what math can do for dance.  

When Wilson first asked me (Wolfson) to the rehearsal studio, I focused on the concept of {\em symmetry}, i.e.
\begin{enumerate}
	\item a collection of transformations, and
	\item a feature left unchanged by those transformations.
\end{enumerate}
In fractals, the transformation takes the form of scaling (some parameter), and the symmetry is then manifested as repetition or similarity across different scales. While symmetry has been an organizing principle in math for centuries, it’s only in the last 50 years that mathematicians have seriously studied fractal symmetry.  Many cultures, especially in Africa and the African diaspora, have developed rich knowledge of how to generate and use fractal symmetry to construct objects, rhythms, movement structures and meaning. However, Africa is not where western mathematics has historically looked for such knowledge and innovation.

It was surprising to me (Wolfson) to find such examples of fractal symmetry in the source materials I (Wilson) shared.  With the dancers of Fist and Heel Performance Group, we explored how these were constructed, how they could be played with, and how they could be incorporated into and expand existing choreographic knowledge and possibilities.  This appeared to be an exciting new application of mathematics!

\paragraph*{II. Braids.} In 2018, I (Wilson) asked Wolfson back into the rehearsal studio to help decipher patterns and structures in reconstructed Shaker (a Utopian millenarian American religious sect) dances that the company was incorporating into a new performance. One specific movement sequence was giving the dancers difficulty, and when we tried to track the movements based on the paths of individual dancers, we kept getting confused and tangled up.  Wolfson realized that if you tracked pairs of dancers, then each pair traced ``elementary braids’’ and the movement of all four dancers could be represented as a ``4-stranded braid.''   This was particularly striking because Shaker iconography and theology was heavily imbued with braiding, both as fiber art and as emblem of communal meanings. Wolfson explained his experience to the dancers as follows,
\begin{quote}
Braids provide a key mathematical tool for describing how 2D and 3D spaces are put together.  They arise in multiple mathematical contexts, from motion of particles on a surface (e.g. molecules in a solution, robots on a factory floor, satellites orbiting the earth), to the solutions of families of polynomial equations. In my research, I am trying to show that certain equations don’t admit simple solutions because the associated braids are too complicated. Any time you find a new context where braids arise, you want to understand it.  I felt really excited to realize there were braids in Shaker dances, and pretty complicated ones at that. The braids were the defining structural element of the dance. Someone clearly had knowledge of braids here and they were using it in a context that for me was totally unexpected.  What else did they know?  Maybe there is knowledge there that you’ve missed, or structure you can use to understand another application.
\end{quote}

\paragraph*{III. Embodied Geometries} Our collaboration has deepened our ability to ask questions in our individual fields.  Currently our thinking has shifted to what can dance do for math. 

Thurston’s revolutionary work on 3-dimensional geometry was one of the major developments of 20th century mathematics. Thurston advocated approaching geometry from the basic question: What is it like to live in, really live in, the abstract spaces we are studying?  While Thurston’s call has driven fundamental advances over the last forty years, it has largely been interpreted as an unembodied perspective, e.g. the motion of point particles or simple bodies like tinker toys. Can choreographic knowledge push (advance) this further?

One possible source for this further development might be around the concepts of {\em choreographic time}, {\em choreographic space} and {\em choreographic movement}. These provide a theoretical core around which I (Wilson) build dance compositions. I have come to work with choreographic time defined as how long it takes ``something'' to happen. Choreographic space is where ``the thing'' happens, which can include {\em far} (beyond the sphere of touch) and {\em near} (from the tips of your fingers to the inner layer of skin), as well as {\em inner} (one’s internal sense of the body). Choreographic movement is what the body does. In the studio, these act as tools to identify which element the maker/choreographer is trying to manipulate or draw the audience’s attention to (and away from).

I (Wolfson) want to understand how we might make mathematical sense of this embodied and perceptual geometry!   I’m still working to understand how these choreographic concepts function for dancemakers and what mathematical uses they might have.  As Wilson explained to me, in his movement lineage-of-understanding, choreographic time/space/movement cannot be separated, just as we cannot reliably separate space and time in relativistic physics. Moreover, what matters to dancemakers is not that these choreographic concepts exist in a platonic sense, but rather that they serve as powerful compositional devices to generate dance. Yet, despite this, I wonder if these choreographic concepts do point to some underlying ``absolute'' phenomenon (akin to the fixed physical spacetime background).  If so, then how such an absolute structures, interacts with and recedes from people’s attention is very much part of the theoretical framework we’re trying to uncover.  

Another vein to explore may be the tension between what is choreographically critical in the creation of a piece and what knowledge is actually available to an audience. As Wilson explained to me,
\begin{quote}
	In the performing of the dance, we are altering people’s perception. What tends to happen, similar to [our first experiments] with fractals and dance: when I tried to create dance using fractal devices, the fractals seemed to be invisible.  Perception/attention seems counter-intuitive.  There is often a divide between how reality ``is'' structured, and how it feels and appears to the viewer.
\end{quote}
Mathematically, how might we encompass and understand the multiple positions and relationships in this divide?

\section*{Ongoing collaborations}
We now share ideas and talk with each other regularly in the course of our professional lives. Both of us have vested interests in discovering the many connections and potentials for math and dance collaborations. For example, I (Wolfson) have collaborated with \href{https://www.math.uci.edu/~apantano/index.html}{Alessandra Pantano} (UC Irvine, Dept. of Math) and \href{https://drama.arts.uci.edu/faculty/tara-rodman}{Tara Rodman} (UC Irvine, Dept. of Drama) to design and implement \href{https://drive.google.com/file/d/1MS7GujCzuUVl9C3I1UIyt8eYAITOs5Ni/view}{lessons for middle school students on symmetry and dance} through UCI's \href{https://sites.ps.uci.edu/mathceo/}{MathCEO} program.  I (Wilson) am currently advocating for the value of what I call “a choreographic way of seeing”: sharing the ``How'' of training the eye to perceive what is actually there.  Last winter, I co-taught  a \href{https://africana.sas.upenn.edu/node/7963}{graduate seminar} at the University of Pennsylvania with \href{https://anthropology.sas.upenn.edu/people/deborah-a-thomas}{Deborah Thomas} (UPenn, Dept. of Anthropology) on  kinesthetic anthropology, focusing on my particular understanding of Choreographic Time, Space and Movement.  Working with scholars and scientists has clarified the value of this training and ability not only for choreographers, but also with potential application across many disciplines.

I (Wilson) continue to expand the use of math in my choreographic work. I have used my collaboration with Wolfson to build a case for dance-math exploration as an essential part of the Fist and Heel Institute, an ``institute without walls'' built to house, articulate and spread the knowledge and movement practices I have researched, developed and incorporated into live performance with Fist and Heel Performance Group. 

Excitingly, our collaboration is not the only one we’re aware of between dance-makers/thinkers and mathematical scientists.  Other examples include: \href{https://campuspress.yale.edu/emilycoates/}{Emily Coates'} (Yale University, Dept. of Drama) and \href{https://demerslab.yale.edu/}{Sarah Demers'} (Yale University, Dept. of Physics) collaboration exploring the intersection of physics and dance, including a book of the same name \cite{CD}; Raleigh, NC's \href{https://www.blackboxdancetheatre.org/}{Black Box Dance Theater's} and \href{https://sites.google.com/ncsu.edu/tlid/home}{Tye Lidman's} (NC State, Dept. of Math)  collaboration has led to the concert dance performance \href{https://www.blackboxdancetheatre.org/applied-dance-and-inspired-mathematics}{Inspired Mathematics}; and scientist/choreographer/dancer \href{http://catiecuan.com/sgwdv7ufxk9lwgatwvwmell5ha6538}{Catie Cuan} (Stanford, Mechanical Engineering PhD Candidate) researches and creates at the intersection of dance and robotics. All of these feature rich interaction between dance-makers/doers and the mathematical sciences, and each has developed distinct strategies, practices and outcomes. 

Just before the pandemic, we were finishing a proposal to bring together researcher-practitioners operating at the interface of dance and math for a first-ever workshop to highlight, compare and foment conversations and collaborations across these disciplines.  The pandemic has put these (and many other) plans on hold.  We continue to engage and explore how math and dance can inform each other and drive mutual creation and discovery.  We are excited to see new collaborations, experiments and creations as they bubble up.

\paragraph*{Acknowledgements} We thank Raja Feather Kelly, Susan Manning, Tara Rodman and the anonymous reviewer for helpful comments, and Rhetta Aleong for many helpful comments, suggestions and edits on an early draft.

\end{document}